\documentclass[12pt,oneside]{amsart}

\usepackage{amssymb}
\usepackage[all]{xypic}
\usepackage[latin1]{inputenc}
\usepackage{verbatim}

\def\Z{\mathbb{Z}}
\def\1{^{-1}}
\def\iff{\Longleftrightarrow}

\def\vf{\varphi}

\def\G{\Gamma}

\newtheorem{theorem}{Theorem}
\newtheorem{cor}{Corollaire}
\newtheorem{lem}{Lemme}

\title{Finite extension of group with infinite conjugacy classes}

\begin{document}
\maketitle
\begin{center}
{\sc Jean-Philippe PR\' EAUX}\footnote[1]{Centre de recherche de l'Ecole de l'air, Ecole de l'air, F-13661 Salon de
Provence air}\ \footnote[2]{Centre de Math\'ematiques et d'informatique, Universit\'e de Provence, 39 rue
F.Joliot-Curie, F-13453 marseille
cedex 13\\
\indent {\it E-mail :} \ preaux@cmi.univ-mrs.fr}
\end{center}

\begin{abstract}
We give a characterization of the group property of being with infinite conjugacy classes (or {\it icc}, {\it i.e.}
$\not= 1$ and of which all
conjugacy classes except 1 are  infinite) for finite extensions of group.\\
\end{abstract}

\section*{Introduction}

A group is said to be with {\sl infinite conjugacy classes} (or {\sl icc}) if it is non trivial, and if all its
conjugacy classes except $\{ 1\}$ are infinite. This property is motivated by the theory of Von Neumann algebra, since
for any group $\G$, a necessary and sufficient condition for its Von Neumann algebra to be a type $II-1$ factor is that
$\G$ be icc (cf. \cite{roiv}).

The property of being icc has been characterized in several classes of groups : 3-manifolds and $PD(3)$ groups in
\cite{aogf3v}, groups acting on Bass-Serre trees in \cite{ydc} and wreath product of groups in \cite{wpicc}. We will
focus here on groups defined by a finite extension or containing a proper finite index subgroup.

Towards this direction particular results are already known. In \cite{aogf3v} have been proved the following results
:\smallskip\\
{\sl -- Let $G$ be defined by a finite exact sequence :
$$1\longrightarrow K\longrightarrow G{\longrightarrow}\Z_2\longrightarrow 1$$
Then $G$ is icc if and only if $K$ is icc and $G\not= K\times \Z_2$. It is easily seen that the condition can be
rephrased as $K$ icc and the natural homomorphism $Q\longrightarrow Out(K)$ is
injective. }\smallskip\\
{\sl -- Let $G$ be a finite extension of $K$, with $G$
torsion-free. Then $G$ is icc if and only if $K$ is icc.}\\

We give  generalizations of these results. We propose the following characterization of finite extensions of groups
with infinite conjugacy classes. A refined and more general version will also be
given in \S 2.\smallskip\\
 \noindent{\bf Main theorem.} \indent {\sl Let $G$ be a group defined by a finite  extension  :
$$1\longrightarrow K\longrightarrow G{\longrightarrow}Q\ \text{(finite)}\longrightarrow 1$$
Then $G$ is icc if and only if $K$ is icc and the natural homomorphism $\theta : Q \longrightarrow Out(K)$ is
injective.
}\smallskip\\

\section{Proof of the main result}

Let us first fix some notations : if $G$ is a group and $u,v$ are element of $G$, then $u^v$ is the element of $G$
defined by $u^v=v^{-1}uv$. If $H$ is a subgroup of $G$, then $u^H=\{u^v\ |\ v\in H\}$ ; in particular $u^G$ denote the
conjugacy class of $u$ in $G$. Note that the cardinal of $u^G$ equals the index of $Z_G(u)$ in $G$ so that $G\not=1$ is
icc if and only if for any $u\not=1 \in G$, $Z_G(u)$ has an infinite index in $G$.

Let $\pi : G\longrightarrow Aut(K)$ be the homomorphism defined by $\pi(g)(k)=k^g$ for any $k\in K$. It makes the
following diagram commute :
\[\xymatrix{
1\ar[r]&K\ar[r]&G\ar[r] \ar_{\pi}[d] & Q\ar[r]\ar^{\theta}[d] & 1\\
1\ar[r]&Inn(K)\ar[r] &Aut(K)\ar[r] & Out(K)\ar[r]&1
 }
\]
%

\noindent{\sl Proof of the main theorem.} We proceed in several steps.\smallskip\\
{\sl Step 1 : $G$ icc $\implies$ $K$ icc.}\\
 Suppose $K$ is not icc, that is  there exists
$u\not=1$ in $K$ such that $Z_K(u)$ has a finite index in $K$. Since $K$ has finite index in $G$, $Z_K(u)$ has a finite
index in $G$, and hence $Z_G(u)$ which contains $Z_K(u)$ also has a finite index in $G$, so that $G$ is not
icc.\smallskip\\
{\sl Step 2 : $G$ icc $\implies$ $\theta$ injective.}\\
Suppose that $\theta : Q\longrightarrow Out(K)$ is non injective, that is there exists $g\in G\setminus K$ and $h\in K$
such that $\forall k\in K$, $g^{-1}kg=h^{-1}k h$. Then $K$ is contained in $Z_G(gh^{-1})$ and hence $Z_G(gh^{-1})$ has
a
finite index in $G$, so that $G$ is not icc.\smallskip\\
{\sl Step 3 : $G$ not icc $\implies$  $K$ not icc or $\theta$ non injective.}\\
We can suppose $G\not=1$ because otherwise $K=1$ is not icc. Let $g\not=1$ in $G$ such that $Z_G(g)$ has a finite index
in $G$, it must exist since $G\not=1$ is not icc. If $g$ lies in $K$ then $Z_K(g)=Z_G(g)\cap K$ has a finite index in
$K$, so that $K$ is not icc. So suppose in the following that $g$ lies in $G\setminus K$.  Let $H=Z_G(g)\cap K$, it is
a finite index subgroup of $K$, and let $k_0,k_1,\ldots ,k_n$ be a set of representative of $K\mod H$. Let
$N=\bigcap_{i=1}^{n}k_i^{-1}Hk_i$, the normalized of $H$ in $K$ ; $N$ is a finite index normal subgroup of $K$ lying in
$Z_G(g)$.  Consider the centralizer $Z_K(N)$ of $N$ in $K$. If $Z_K(N)\not=1$ let $u$ be a non trivial element of
$Z_K(N)$. Since $u^N$ is a singleton and $N$ has a finite index in $K$ it follows that $u^K$ is finite so that $K$ is
not icc. If $Z_K(N)=1$, once have been noted that  $\pi(g)$  restricted to $N$ is the identity, the following lemma
applies to show that $\pi(g)$ is the identity of $K$, so that $\theta$ is non injective.\hfill $\square$

\begin{lem} Let $K$ be a group, $N$ a normal subgroup of $K$ and $\vf\in Aut(K)$ which is the identity once restricted
to $N$. If $Z_K(N)=1$ then $\vf$ is the identity on $K$.
\end{lem}

\noindent{\sl Proof.}  Let $k\in K$, then for any $h\in N$,
$$khk^{-1}=\vf(khk^{-1})=\vf(k)h\vf(k^{-1})$$
so that $k^{-1}\vf(k)\in Z_K(N)$. Hence if $Z_K(N)=1$, for any $k\in K$, $\vf(k)=k$.\hfill$\square$

\section{Finite index subgroups in icc groups}

Let $H$ be a subgroup of the group $G$ ; the {\sl normalized} $N(H)$ of $H$ is defined to be $N(H)=\bigcap_{g\in G}
H^g$. It is a normal subgroup both in $H$ and $G$, and $N(H)=H$ exactly when $H$ is normal in $G$. If $H$ has a finite
index in $G$, let $g_1,g_2,\ldots , g_n$ be a finite set of representatives of $G\mod H$, then
$N(H)=\bigcap_{i=1}^nH^{g_i}$, and $N(H)$ has a finite index both in $H$ and in $G$. Pay attention that the normalizeD
has nothing to do with the normalizeR.

Denote by $\pi$ the homomorphism $\pi : G\longrightarrow Aut(N(H))$ defined by $\forall g\in G, k\in N(H)$,
$\pi(g)(k)=k^g$.

\begin{theorem} Let $G$ be a group, $H$ a finite index subgroup of $G$ and $N(H)$ the normalized of $H$.
 Then $G$ is not icc if and only if at least one of the conditions (i), (ii) or (iii) is satisfied : \smallskip\\
(i) \indent $H$ is not icc\smallskip\\
(ii) \indent $\exists\ g\not=1\in G\setminus H$ with finite order such that $\pi(g)$ is the identity,\smallskip\\
(iii)\indent $\exists\ g\not=1\in G\setminus H$ with finite order such that $\pi(g)$ is inner.\smallskip\\
\indent Moreover on the one hand (ii) $\implies $ (iii) and on the other (iii) $\implies$ (i) or (ii), so that :\\
\centerline{$G$ is not icc $\iff$ (i) or (ii) $\iff$ (i) or (iii).}
\end{theorem}

\noindent{\sl Proof of theorem 1.} We proceed in several steps :\smallskip\\
{\sl Step 1. $(i)$ $\implies$ $G$ not icc.}\\ Since $H$ has a finite index in $G$, if $H$ is not icc then $G$ is
clearly also not icc.\smallskip\\
{\sl Step 2. (ii) $\implies$ $G$ not icc.}\\ If $\exists\ g\not=1\in G\setminus H$ with finite order $n$ such that
$\pi(g)$ is the identity, then $G$ contains $N(H)\times Z_n$ which is obviously not icc and has a finite index in $G$,
so that $G$ is not icc.\smallskip\\
{\sl Step 3. (iii) $\implies$ (i) or (ii).}\\ Let $g\not=1\in G\setminus H$ be a finite order element such that
$\exists k\in N(H)$, $\forall h\in H$, $\pi(g)(h)=h^k$. Then $\pi(gk^{-1})$ is the identity on $N(H)$. Let $n\in \Z$ be
the order of $g$, then $(gk^{-1})^n$ lies in the center of $N(H)$. Hence either $N(H)$ has a non trivial center and
condition $(i)$ follows, or $(gk^{-1})^n=1$ so that condition $(ii)$ is satisfied.\smallskip\\
{\sl Step 4. $G$ not icc $\implies$ (i) or (ii).}\\
Suppose $G$ is not icc ; $G$ is a finite extension of $N(H)$. Either $H$ is not icc or $\exists\ g\in G\setminus H$
with $g^G$ finite. For some $n>1$, $g^n\in N(H)$ ; if $g^n\not= 1$ then $N(H)$ is not icc and it follows that condition
$(i)$ is satisfied. So suppose in the following that $g^n=1$. The same argument as in the step 3 of the proof of the
main theorem --with $N(H)$ instead of $K$-- shows that either $N(H)$ is not icc, so that condition $(i)$ is satisfied,
or $\pi(g)$ is the identity, so that condition $(ii)$ is satisfied.\hfill$\square$\\

%
%
\begin{cor}
$G$ is icc if and only if $H$ is icc and $G\not\supset N(H)\times \Z_n$.
\end{cor}
%
%
\begin{cor}
If $G\setminus H$ contains no torsion element (in particular when $G$ is torsion-free) then $G$ is icc if and only if
$H$ is icc.
\end{cor}
\newpage


\begin{thebibliography}{mot}

\bibitem[Co]{ydc}
Y.\textsc{de Cornulier}, \emph{Infinite conjugacy classes in groups acting on trees}, preprint (2005).


\bibitem[HP]{aogf3v}
P.\textsc{de la Harpe} et J.-P.\textsc{Préaux}, \emph{Groupes fondamentaux des variétés de dimension 3 et algèbres
d'opérateurs}, preprint arXiv:math.GR/0509449 v1 (2005).

\bibitem[Pr]{wpicc}
J.-P.\textsc{Préaux}, \emph{Wreath product of groups wit infinite conjugacy classes}, preprint arXiv:math.GR/0612685
(2006).

\bibitem[ROIV]{roiv}
F.J.\textsc{Murray} et J.\textsc{von Neumann}, \emph{On rings of operators}, IV, Annals of Math. {\bf 44} (1943),
716--808.

\end{thebibliography}
 \end{document}